\documentclass[11pt]{article}
\usepackage{amsmath}
\usepackage{amssymb}

\renewcommand{\Bbb}{\mathbb}
\newcommand{\B}{{\cal  B}}
\newcommand{\C}{{\Bbb  C}}
\newcommand{\D}{{\Bbb D}}

\newcommand{\R}{{\Bbb  R}}

\renewcommand{\P}{{\Bbb  P}}

\newcommand{\T}{{\Bbb  T}}
\newcommand{\A}{{\cal A}}

\newcommand{\LL}{{\cal L}}

\newcommand{\OO}{{\cal O}}

\newcommand{\PSH}{{\operatorname{{\cal PSH}}}}

\renewcommand{\phi}{\varphi}
\renewcommand{\epsilon}{\varepsilon}

\newcommand{\eq}{\begin{equation}}
\newcommand{\ee}{\end{equation}}

\def\squarebox#1{\hbox to #1{\hfill\vbox to #1{\vfill}}}

\newtheorem{theorem+}           {Theorem}      [section]
\newtheorem{definition+}  [theorem+]  {Definition}
\newtheorem{lemma+}  [theorem+]  {Lemma}
\newtheorem{corollary+}  [theorem+]  {Corollary}
\newtheorem{proposition+}  [theorem+]  {Proposition}
\newtheorem{example+}  [theorem+]  {Example}
\newtheorem{question+}  [theorem+]  {Question}

\newenvironment{proof}{\medbreak\noindent{\it Proof:}\rm}{\hfill$\square$\rm}

\title{\Large \bf DISC FORMULAS FOR THE WEIGHTED SICIAK-ZAHARIUTA
EXTREMAL FUNCTION\footnote{\it Dedicated to J\'ozef Siciak on the occasion of his 75th
birthday}
}
\author{\large\bf  Benedikt Steinar Magn\'usson  and  Ragnar Sigurdsson }

\date{December 29, 2006}

\begin{document}
\maketitle

\begin{abstract}  \noindent
We prove  a disc formula for the weighted Siciak-Zahariuta extremal
function  $V_{X,q}$ for an upper semicontinuous function $q$ 
on an open connected subset $X$ in $\C^n$.  This function is also known as the
weighted Green function with logaritmic pole at infinity
and weighted global extremal function.


\end{abstract}

\bigskip
{\bf Introduction.} \ 
If  $X$ is a subset of $\C^n$ and $q:X\to
\overline\R=[-\infty,+\infty]$ is a  function, then the   
{\it weighted Siciak-Zahariuta
extremal function} $V_{X,q}$ {\it with respect to } $q$ is defined as
$$
V_{X,q}=\sup\{u\in \LL\,;\, u\leq q \text{ on } X\},
$$
where $\LL$ denotes the Lelong class, which consists of all
plurisubharmonic functions $u$ on $\C^n$ of minimal growth, i.e.,
functions $u$ satisfying $u(z)\leq \log^+\|z\|+c_u$, $z\in \C^n$ for
some constant $c_u$.  The {\it Siciak-Zahariuta extremal function}
$V_X$ corresponds to the case $q=0$.  The functions $V_X$ and
$V_{X,q}$  were first introduced by Siciak in the fundamental paper
\cite{Sic} where he proved his celebrated approximation theorem in
several complex variables.  The theorem states that for every compact
subset $X$ of $\C^n$, such that $V_X$ is continuous, a holomorphic
function $f$ on some neighbourhood of $X$ can be approximated
uniformly on $X$ by polynomials $P_\nu$ of degree less than or equal
to $\nu$ in such a
way that
$$
\limsup_{\nu\to \infty}\big(\sup_{z\in
E}|f(z)-P_\nu(z)|\big)^{1/\nu}=\varrho<1
$$
if and only if $f$ has a holomorphic extension to the sublevel
set $\{z\in \C^n\,;\, V_X(z)<-\log \varrho\}$.

The purpose of this paper is to extend the methods of
L\'arusson and Sigurdsson \cite{LarSig} in order to prove 
disc envelope formulas for  $V_{X,q}$.  Our main result is the
following

\medskip\noindent
{\bf Theorem 1.} \ {\it  Let $X$ be an open connected subset of $\C^n$ and  $q$ be an
upper semicontinuous  function on $X$.  Then for every $z\in \C^n$
$$
V_{X,q}(z)=\inf\big\{-\sum_{a\in f^{-1}(H_\infty)}
\log|a|+\int_\T q\circ f\, d\sigma \,;\,  f\in \OO(\overline
\D,\P^n), f(\T)\subset
X, f(0)=z\big\}. 
$$
Here $\P^n$ is the complex projective space viewed in the usual
way as the union of the affine space $\C^n$ and the hyperplane at
infinity $H_\infty$, $\D$ and $\T$ are the open unit disc and  the unit
circle in $\C$, and $\sigma$ is the normalized arc length measure on 
$\T$.
}

\medskip
Our approach in the paper is the following.  Based on the observation,
see Guedj and Zeriahi \cite{GueZer},
that a function $u$ is in the Lelong class if and only if 
$(z_0,\dots,z_n)\mapsto u(z_1/z_0,\dots,z_n/z_0)+\log|z_0|$ extends
as a plurisubharmonic function from $\C^{n+1}\setminus\{z_0=0\}$ to
$\C^{n+1}\setminus\{0\}$, we derive a fundamental inequality
$u(z)\leq J_q(f)$, for any closed analytic disc mapping 
the origin to $z$ and the unit circle into $X$.  This inequality 
defines a disc functional $J_q$ associated to $q$.  Then we define the 
{\it good } sets of analytic discs with respect to $q$ and observe that Poletsky's
theorem implies a disc formula for $V_{X,q}$.  From this formula 
we prove that $V_{X,q}$ is the envelope of $J_q$ with respect to the
class of all closed analytic discs mapping the unit circle into
$X$.  This result gives the theorem above.

\bigskip 
{\bf Notation and some basic results.} \ 
An analytic disc in a manifold $Y$ is a holomorphic map
$f:\D\to Y$ from  the unit disc $\D$ in $\C$ into $Y$.
We denote the set of all analytic discs in $Y$ by $\OO(\D,Y)$.
A disc functional on $Y$ is a map $H:\A\to \overline \R$ defined on some 
subset $\A$ of $\OO(\D,Y)$ with values in the extended real
line  $\overline \R=[-\infty,+\infty]$.  The envelope $E_{\B}H:Y\to \overline \R$ 
of $H$ with respect to the subclass $\B$ of $\A$ is defined by
$$
E_{\B}H(x)=\inf\{H(f)\,;\, f\in \B, f(0)=x\}, \qquad x\in Y.
$$
We let $\A_Y$ denote the set of all closed analytic discs in $Y$,
i.e., analytic discs that extend to holomorphic maps in some
neighbourhood of the closed unit disc $\overline \D$, and 
for a subset $S$ of $Y$ we let $\A_Y^S$ denote the set of
all discs in $\A_Y$ which map the unit circle $\T$ into 
$S$.

We let $\P^n$ denote the complex projective space with the natural 
projection $\pi:\C^{n+1}\setminus \{0\}\to \P^n$,
$(z_0,\dots,z_n)\mapsto [z_0:\cdots:z_n]$ and we identify
$\C^n$ with the subspace of $\P^n$ consisting of all 
$[z_0:\cdots:z_n]$ with $z_0\neq 0$.  The hyperplane at infinity
$H_\infty$ in $\P^n$ is the projection of $Z_0\setminus \{0\}$
where $Z_0$ is the hyperplane in $\C^{n+1}$ defined by the equation 
$z_0=0$.  

It is an easy observation that a function $u\in \PSH(\C^n)$ is in the
Lelong class $\LL$ if and only if the function  
\begin{equation}\label{eq:1}
\tilde z=(z_0,\dots,z_n)\mapsto u\circ \pi(\tilde z)+\log|z_0|
=u(z_1/z_0,\dots,z_n/z_0)+\log|z_0|
\end{equation}
extends as a plurisubharmonic function from $\C^{n+1}\setminus Z_0$ 
to $\C^{n+1}\setminus \{0\}$.  If we denote this extension by $v$, 
take $f=[f_0:\cdots:f_n]\in \A_{\P^n}$ with $f(0)=z\in \C^n$,
$f(\T)\subset \C^n$,
and set $\tilde f=(f_0,\dots,f_n)\in \A_{\C^{n+1}\setminus \{0\}}$,
then by subharmonicity of  $v\circ \tilde f$ we get
\begin{equation}\label{eq:2}
u(z)+\log|f_0(0)|=v\circ \tilde f(0) \leq \int_\T v\circ \tilde f\, d\sigma
=\int_\T u\circ f\, d\sigma+\int_\T\log|f_0| \, d\sigma.
\end{equation}
Since  $f(\T)\subset \C^n$, the set $f(\D)$ has finitely many
intersections with $H_\infty$, which means that $f_0$ has finitely
many zeros in $\D$.  We write
$$
f_0(\zeta)=\prod_{a\in f^{-1}(H_\infty)} \bigg(\dfrac{\zeta-a}{1-\bar
a\zeta}\bigg)^{m_{f_0}(a)} \, g_0(\zeta),
$$
where $m_{f_0}(a)$ denotes the multiplicity of $a$ as a zero of
$f_0$ and $g_0$ is holomorphic and without zeros in some neighbourhood 
of $\overline \D$.  We have 
\begin{equation}\label{eq:3}
\log|f_0(0)|=\sum_{a\in f^{-1}(H_\infty)}
m_{f_0}(a)\log|a|+\log|g_0(0)|,
\end{equation}
and since the product has modulus  $1$ on $\T$ and $\log|g_0|$ is
harmonic in some neighbourhood of $\overline \D$, we have
\begin{equation} \label{eq:4}
\int_\T\log|f_0|\, d\sigma=\int_\T\log|g_0|\, d\sigma=\log|g_0(0)|.
\end{equation}
By combining (\ref{eq:3}) and (\ref{eq:4}) with (\ref{eq:2}) we
arrive at the inequality
\begin{equation}\label{eq:5}
u(z)\leq -\sum_{a\in f^{-1}(H_\infty)} m_{f_0}(a)\log|a|
+\int_\T u\circ f\, d\sigma.
\end{equation}
As in \cite{LarSig} we define the disc functional 
$$
J:\OO(\D,\P^n)\to \overline \R_+=[0,+\infty], \qquad
J(f)=-\sum_{a\in f^{-1}(H_\infty)} m_{f_0}(a)\log|a|,
$$
where we take $J(f)=0$ if $f^{-1}(H_\infty)=\varnothing$.
If $q$ is Borel measurable, then we add a mean value term to 
$J$ and define $J_q$ by
$$
J_q:\OO(\D,\P^n)\cap C(\overline \D,\P^n) \to \overline \R, \qquad
J_q(f)=J(f)+\int_{\T\cap f^{-1}(X)} q\circ f\, d\sigma.
$$
If $f^{-1}(H_\infty)$ is an infinite set the sum is taken 
as the infimum over all finite subsets, which is well defined
since the terms are all negative.  In the case when $J(f)=+\infty$
and the integral is $-\infty$ we define $J_q(f)=+\infty$.
If $f(\T)\subset X$, then the sum is finite.
For the constant disc $k_x$, $\overline \D\ni \zeta\mapsto x\in X$, 
we have $J(k_x)=0$, and hence $J_q(k_x)=q(x)$.

The inequality (\ref{eq:5}) implies that for every $u\in \LL$ with 
$u\leq q$ on $X$ and every  $f\in \A_{\P^n}$
with $f(0)=z$ we have
$$
u(z)\leq J_q(f)+\int_{\T\setminus f^{-1}(X)} u\circ f\, d\sigma.
$$
If  $f(\T)\subset X$, then the second term in the right hand side
vanishes.    If we take the supremum over all $u\in \LL$ with
$u\leq q$ on $X$ in the left hand side and the infimum over all 
$f\in \B$ for some subclass  $\B\subseteq \A_{\P^n}^X$ in the right 
hand side, then we arrive at the inequality
$$
V_{X,q}(z)\leq E_{\A_{\P^n}^X}J_q(z)\leq E_\B J_q(z),  \qquad z\in
\C^n.
$$
We will prove that the first inequality
is actually an equality:

\bigskip
{\bf Theorem 2.} {\it  Let $X$ be an open connected subset of $\C^n$
and   $q:X\to
\R\cup \{-\infty\}$ be an upper semicontinuous function.  Then
$V_{X,q}=E_{\A_{\P^n}^X}J_q$, i.e., for every $z\in \C^n$ we have
$$
V_{X,q}(z)=\inf\{
-\sum_{a\in f^{-1}(H_\infty)} m_{f_0}(a)\log|a|
+\int_{\T} q\circ f\, d\sigma
  \,;\, f\in \A_{\P^n}, f(\T)\subset X, f(0)=z\}.
$$
}

\bigskip
Observe that the formula in Theorem 1 is 
the same as this one  except for the multiplicities.  
In order to show that Theorem 1 follows from Theorem 2,
we first observe that the upper semicontinuity of $q$ implies that
for every $\varepsilon>0$ and every $f\in \A_{\P^n}^X$
there exists a continuous function $\tilde q\geq q$ on $X$
such that $\int_\T\tilde q\circ f\, d\sigma 
< \int_\T q\circ f\, d\sigma +\varepsilon$.  By Proposition
1 in \cite{LarSig}, every $f\in \A_{\P^n}$ can be approximated
uniformly on $\overline \D$ by $g\in \A_{\P^n}$, such that
all the zeros of $g_0$ are simple, $g(0)=f(0)$, and $J(g)=J(f)$.
Since $\tilde q$ is continuous we can choose $g$ such that
$\int_\T\tilde q\circ g\, d\sigma <\int_\T\tilde q\circ f\, d\sigma
+\varepsilon$.  This gives that $J_q(g)\leq J_{\tilde
q}(g)<J_q(f)+2\varepsilon$ and we conclude that the infima in
Theorem 1 and 2 are equal.

\bigskip
{\bf Good sets of analytic discs.} \ 
We modify the definition from \cite{LarSig} of {\it good sets} of analytic discs by 
saying that a subset $\B$ of $\A_{\P^n}$ is 
{\it good with respect to the function} $q$ if:

\begin{description}
\item{(1)} $f(\T)\subset X$ for every $f\in\B$,
\item{(2)} for every $z\in\C^n$, there is a disc in $\B$ with centre $z$,
\item{(3)} for every $x\in X$, the constant disc at $x$ is in $\B$, and
\item{(4)} the envelope $E_\B J_q$ is upper semicontinuous on $\C^n$ and has
minimal growth, that is, $E_\B J_q-\log^+\|\cdot\|$ is bounded above on
$\C^n$.
\end{description}

\noindent
The condition (1) implies that $u(z)\leq J_q(f)$ for every $u\in \LL$
with $u\leq q$ and $f\in\B$ with $f(0)=z$, (2) implies that 
$E_\B J_q(z)<+\infty$ for every $z\in \C^n$,
(3) implies that $E_\B J_q(x)\leq q(x)$ for all $x\in X$, and 
(4) implies that $V_{X,q}$ is the largest plurisubharmonic function
on $\C^n$ dominated by $E_\B J_q$.    

Poletsky's theorem states that for every upper semicontinuous function
$\psi:Y\to \R\cup \{-\infty\}$ on a complex manifold $Y$ we have 
$$
\sup\{u(x)\,;\, u\in \PSH(Y), u\leq \psi\}
=\inf\{\int_\T \psi\circ h\, d\sigma \,;\, h\in \A_Y, h(0)=x\},
\qquad x\in Y.
$$
See Poletsky \cite{Pol}, L\'arusson and Sigurdsson \cite{LarSig1, LarSig2}, and Rosay \cite{Ros}.
As a consequence we get a disc formula for $V_{X,q}$:

\bigskip
{\bf Theorem 3.} {\it  Let $X$ be an open subset of $\C^n$,  $q:X\to
\overline \R$ be a  Borel measurable function, and $\B$ be a good class
of analytic discs with respect to $q$.  Then
$$
V_{X,q}(z)=\inf\{\int_\T E_\B J_q \circ h \, d\sigma  \,;\, h\in \A_{\C^n}, h(0)=z\},
\qquad z\in \C^n.
$$
}

\bigskip
{\bf The remaining proof.} \  Assume that $q:X\to \R\cup\{-\infty\}$
is upper semicontinuous.
From now on we will choose  $\mathcal B$ as the set of all analytic
disc  in $\P^n$ which are either a constant
disc in $X$ or of the following form
$$
f_{z,w,r} : \zeta \mapsto w + \frac{\| z-w \| + r\zeta}
{r+\| z-w \| \zeta} \frac{r}{\| w-z \|} (z-w),
$$
where $z \in \C^n$, $w \in X\setminus \{z\}$ and $r < \min \{
\|z-w\|, d(w,\partial X)\}$.

Observe that $f_{z,w,r}$ maps $\overline \D$ into the projective line
through $z$ and $w$, $\T$ is mapped to the circle with centre
$w$ and radius $r$, $0$ is mapped to $z$, and $-r/\|z-w\|$ is mapped
into $H_\infty$.  The conditions on $z$, $w$ and $r$ ensure that
$f_{z,w,r}(\T)\subset X$ and we have the formula
\begin{equation}\label{eq:6}
J_q(f_{z,w,r})=\log(\|z-w\|/r)+\int_{\T} q\circ f_{z,w,r}\, d\sigma.
\end{equation}
It is obvious that the conditions (1), (2), and (3) in the definition
of a good set are satisfied.
We have by (3) that $E_\B J_q(x)\leq q(x)$ for all $x\in X$, and since
$q$ is upper semicontinuous, this implies that $E_\B J_q$ upper
bounded on every compact subset of $X$.  If we fix $w\in X$ and
$r<d(w,\partial X)$, then it follows from (\ref{eq:6}) that 
$E_\B J_q$ is upper bounded on every compact subset of $\C^n$ and is
of minimal growth.   The upper semicontinuity of $E_\B J_q$ follows
from

\bigskip
{\bf Lemma 1.} \ {\it Assume that $q:X\to \R\cup\{-\infty\}$
is upper semicontinuous.  For every $z_0\in \C^n$ and every 
$\alpha\in \R$ such that 
$E_\B J_q(z_0)<\alpha$ there exist $w_0\in \C^n$, $r_0>0$, and  
a neighbourhood $U$ of $z_0$ such that
$0<r_0<\min\{\|z-w_0\|,d(w_0,\partial X)\}$
and $J_q(f_{z,w_0,r_0})<\alpha$ for all $z\in U$.}

\begin{proof}  Let $f\in \B$ such that $f(0)=z_0$ and $J_q(f)<\alpha$.
If $f$ is of the form $f_{z_0,w_0,r_0}$ for some $w_0\in \C^n$ and
$0<r_0<\min\{d(w_0,\partial X),\|z_0-w_0\|\}$, then we can choose a
continuous function $\tilde q\geq q$ on $X$ such that $J_{\tilde
q}(f_{z_0,w_0,r_0})<\alpha$.  The continuity of $\tilde q$ implies 
that there exists a neighbourhood $U$ of $z_0$ such that 
$r_0<\|z-w_0\|$ and $J_{\tilde q}(f_{z,w_0,r_0})<\alpha$ for all 
$z\in U$. Since $J_q\leq J_{\tilde q}$ the statement holds in this case.

Assume now that $f$ is the constant disc $z_0$.   Then $z_0\in X$ and
$J_q(f)=q(z_0)<\alpha$.  Since $q$ is upper semicontinuous, there
exists $0<\delta<d(z_0,\partial X)$ such that $q(z)<\alpha$ for all 
$z\in B(z_0,\delta)$, the ball with center $z_0$ and radius $\delta$.
Then for every $z$ and $w$ in $B(z_0,\tfrac 12 \delta)$ and 
$0<r<\min\{\|z-w\|,\tfrac 12 \delta\}$ we have $\int_\T 
q\circ f_{z,w,r}\, d\sigma <\alpha$.
Now choose $w_0\in B(z_0,\tfrac 12\delta)$ and $0<r_0<
\min\{\|z_0-w_0\|,\tfrac 12\delta\}$ such that 
$J_q(f_{z_0,w_0,r_0})=\log(\|z_0-w_0\|/r_0)
+\int_{\T} q\circ f_{z_0,w_0,r_0}\, d\sigma<\alpha$.  
The statement now follows as in the first part of
the proof.
\end{proof}    

\bigskip
If $q:X\to \R\cup\{-\infty\}$ is 
upper semicontinuous and $q_j:X\to \R$ is a decreasing sequence 
of continuous functions converging to $q$,
then it is obvious that $V_{X,q_j}\searrow V_{X,q}$. It also 
immediately follows that  $J_{q_j}(f)\searrow J_q(f)$ for 
every $f\in \A_{\P^n}^X$ and as a consequence we get 
$E_{\A_{\P^n}^X}J_{q_j}\searrow  E_{\A_{\P^n}^X}J_{q}$.
This shows that for the proof of Theorem 2 we may assume that
$q$ is  continuous.  

In the previous section we have seen that  $V_{X,q}\leq
E_{\A_{\P^n}^X}J_q$  and that $V_{X,q}$ is the largest
plurisubharmonic function on $\C^n$ dominated by $E_{\B}J_q$.
Hence, Theorem~2 is a direct consequence of Theorem~3 and the following 

\bigskip
{\bf Lemma 2.} \ {\it 
Let $X$ be an open connected subset of $\C^n$, $q:X\to \R$ be
continuous, and $\mathcal B$ be as above. For every $h\in \A_{\C^n}$,
every continuous function $v \geq E_{\mathcal B} J_q$ on $\C^n$,
and every $\varepsilon > 0$, there exists $g \in \mathcal
A_{\P^n}^X$ with  $g(0)=h(0)$ and
$$
J_q(g) \leq \int_\T v\circ h \, d\sigma + \varepsilon.
$$
}

The proof is exactly the same as the proof of the Lemma in
\cite{LarSig} with $J_q$ in the role of $J$.  We only have to note 
that if we choose $\varphi:\C^{n+1}\to \R$ with $\varphi(z)=\log|z_0|$,
let $f=[f_0:\cdots:f_n]\in \A_{\P^n}$, and let  
$\tilde f=(f_0,\dots,f_n)\in \A_{\C^{n+1}\setminus\{0\}}$
be a lifting of $f$, then
$$
J_q(f)=\int_\T(\varphi\circ \tilde f+q\circ \pi\circ \tilde f)\,
d\sigma-\varphi(\tilde f(0)),
$$
and that the last part of the proof holds with $\varphi+q\circ \pi$ in 
the role of $\varphi$.

\bigskip
{\small
Science Institute, University of Iceland, Dunhaga 3, IS-107
Reykjavik, Iceland

E-mail: bsm@hi.is and  ragnar@hi.is
}
\end{document}